\documentclass[10pt, a4paper, UKenglish]{article}

\usepackage[UKenglish]{babel}
\usepackage{amssymb, amsmath, textcomp, amsthm}
\usepackage{comment}
\usepackage{bbm,dsfont}
\usepackage{german}
\usepackage[pdftex]{graphicx}
\usepackage[utf8]{inputenc}
\usepackage{authblk}
\numberwithin{equation}{section}
\title{A geometric proof of Bourgain's $L^2$ estimate of maximal operators along analytic vector fields}
\author{Shaoming Guo}
\date{}

\def\R{\mathbb{R}}
\def\N{\mathbb{N}}

\def\Z{\mathbb{Z}}

\def\lesim{\lesssim}

\def\begineq{\begin{equation}}
\def\endeq{\end{equation}}

\theoremstyle{plain}
\newtheorem{thm}{Theorem}[section]

\newtheorem{lem}[thm]{Lemma}

\newtheorem{defi}[thm]{Definition}
\newtheorem{claim}[thm]{Claim}
\newtheorem{rem}[thm]{Remark}

\newtheorem*{openproblem*}{Open Problem}

\begin{document}
\maketitle

\begin{abstract}
Bourgain \cite{Bo}  proved  that the maximal operator associated to an analytic vector field is bounded on $L^2$. In the present paper, we give a geometric proof of Bourgain's result by using the tools developed by Lacey and Li in \cite{LL1} and \cite{LL2}.
\end{abstract}

\let\thefootnote\relax\footnote{Date: \date{\today}; MSC classes:  42B20, 42B25.}

\section{Statement of the main result}

Let $\Omega\subset \R^2$ be a bounded open set and $v:\Omega'\to S^1$ be a unit vector field defined on a neighbourhood $\Omega'$ of the closure $\bar{\Omega}$ of $\Omega$. For a fixed small positive number $\epsilon_0>0$, define the maximal operator associated to the vector field $v$ truncated at the scale $\epsilon_0$ by
\begineq\label{EE1.1}
M_{v, \epsilon_0}f(x):=\sup_{\epsilon<\epsilon_0}\frac{1}{2\epsilon}\int_{-\epsilon}^{\epsilon}|f(x+tv(x))|dt.
\endeq
In \cite{Bo}, Bourgain proved that
\begin{thm}[\cite{Bo}]\label{theorem1.1}
Let $v$ be real analytic on $\Omega'$. Then for $\epsilon_0>0$ chosen small enough, $M_{v, \epsilon_0}$ is bounded on $L^2(\Omega)$.
\end{thm}

\begin{rem}
The $L^p$ bounds $(\forall p>1)$ for \eqref{EE1.1} and its singular integral variant
\begineq
H_{v, \epsilon_0}f(x):=\int_{-\epsilon_0}^{\epsilon_0}f(x+tv(x))\frac{dt}{t}
\endeq
were obtained by Stein and Street \cite{SS} via a very different method. See also \cite{CNSW} for results concerning certain smooth vector fields. 
\end{rem}

To prove Theorem \ref{theorem1.1}, Bourgain reduced the analyticity assumption on the vector field to the following geometric one: for $x\in \Omega$ and $t$ small enough, define the function
\begineq
\omega_x(t)=\left|\det[v(x+t v(x)), v(x)]\right|.
\endeq
We assume that
\begineq\label{AA1.3}
\left| \{t\in [-\epsilon, \epsilon]: \omega_x(t)< \tau \sup_{-\epsilon\le s\le \epsilon}\omega_x(s)\} \right|\le C_0 \tau^{c_0} \epsilon,
\endeq
for all $0<\tau<1, 0<\epsilon\le \epsilon_0$, where $0<c_0, C_0<\infty$ are constants independent of the point $x\in \Omega$.\\

It is shown in \cite{Bo} that Theorem \ref{theorem1.1} can be reduced to the following
\begin{thm}[\cite{Bo}]\label{theorem1.2}
If $v$ is $C^1$ and satisfies the condition \eqref{AA1.3}, then $M_{v, \epsilon_0}$ is bounded on $L^2(\Omega)$.
\end{thm}

Bourgain's proof for Theorem \ref{theorem1.2} is not entirely geometric, especially in the key Lemma 3.28, where he used polar coordinates and applied Schur's lemma to get the desired $L^2$ bounds.\\

The goal of this paper is to give a (relatively) geometric proof of Theorem \ref{theorem1.2}. The idea is to use the time-frequency decomposition initiated by Lacey and Li in the setting of the Hilbert transform along vector fields in \cite{LL1} and \cite{LL2}. This decomposition was further developed by Bateman \cite{Ba2}, Bateman and Thiele \cite{BT}. We refer to \cite{Guo} for the detailed description of the progress that has been made for Hilbert transforms and maximal operators along vector fields.

This paper is free of time-frequency analysis techniques. Other than certain geometric lemmas (in the following Section \ref{section44}), the tools that we will be using are only the elementary ones, for example Plancherel's identity, the Cauchy-Schwartz inequality, Minkowski's inequality and so on.\\

{\bf Notations:} Throughout this paper, we will write $x\ll y$ to mean that $x\le y/10$, $x\lesim y$ to mean that there exists a universal constant $C$ s.t. $x\le C y$, and $x\sim y$ to mean that $x\lesim y$ and $y\lesim x$.  $\mathbbm{1}_E$ will always denote the characteristic function of the set $E$.\\

{\bf Acknowledgement:} The author would like to thank his advisor, Prof. Christoph Thiele, for sharing with him this problem and for many valuable discussions.

\section{Reduction to a smooth cut-off}\label{section2}
By a renormalization, we assume further that $\|v\|_{C^1}\le 1$, and $\Omega=B_{\epsilon_0}(0)$, which is the ball of radius $\epsilon_0\ll 1$ centered at the origin. Moreover, as we are only concerned with the truncated maximal operator, we can w.l.o.g. assume that the vector field is periodic in both horizontal and vertical directions with each periodicity being $3\cdot \epsilon_0$, and that the vector field always points in the two-ended cone which forms an angle less than $\pi/10$ with the horizontal axis. In the following, we will denote this cone by $\Gamma_0$.\\

Choose $\alpha:\R\to \R$ to be a proper smooth bump function such that the support of $\hat{\alpha}$ lies on $[-1, 1]$. For $0<\epsilon\le \epsilon_0$, define
\begineq\label{AA2.1}
A_{\epsilon}f(x):=\int_{\R}f(x+\epsilon t v(x))\alpha(t)dt.
\endeq
It is not difficult to see that to bound the maximal function $M_{v, \epsilon_0}f$, it is essentially equivalent to bound 
\begineq
\sup_{j\in \N, 2^{-j}\le \epsilon_0}\left|A_{2^{-j}}f(x)\right|,
\endeq
which will still be called $M_{v, \epsilon_0}f$. Moreover we will write $A^j$ to stand for $A_{2^{-j}}$ for the sake of simplicity. Hence in the rest of the paper, we will focus on the following operator
\begineq\label{AA2.3}
M_{v, \epsilon_0}f(x):= \sup_{j\in \N, 2^{-j}\le \epsilon_0}|A^j f(x)|,
\endeq
and prove that it is bounded on $L^2$.

%In the end, we just need to observe that to prove \eqref{AA2.3}, it suffices to prove the following spatially localized version:
%\begineq
%\|M_{v, \epsilon_0}f\|_{L^2(B_{\epsilon_0}(0))}\lesim \|f\|_{L^2},
%\endeq
%due to the fact that we are truncating the maximal operator at the scale $\epsilon_0$.

\section{Bourgain's high-low frequency decomposition}
We linearize the maximal operator in \eqref{AA2.3}: Take a measurable function $J:\R^2\to \N$ such that
\begineq
M_{v, \epsilon_0}f(x)\sim \left|A^{J(x)}f(x)\right|.
\endeq

For a point $x\in \Omega=B_{\epsilon_0}(0)$, let $R_{x, j}$ be the rectangle with center $x$, orientation $v(x)$, length $2^{-j}$ in direction $v(x)$ and width 
\begineq
\delta(R_{x, \epsilon})=2^{-j}\cdot \sup_{|t|<2^{-j}}\omega_x(t).
\endeq
Especially we denote
\begineq
\delta(x):=2^{-J(x)}\cdot \sup_{|t|<2^{-J(x)}}\omega_x(t).
\endeq
Choose a measurable function $K:\R^2\to \N$ such that 
\begineq
\delta(x)\sim 2^{-K(x)}, \forall x\in \R^2.
\endeq
Do a Littlewood-Paley decomposition for the function $f$, and write
\begineq
f=\sum_{k\in \Z}P_k f.
\endeq
This turns our linearised maximal operator into
\begineq
\sum_{k\in \Z}A^{J(x)}P_k f(x).
\endeq
Bourgain's idea is to split the function $f$ into two parts, the high frequency part and the low frequency part, in the following way:  
\begineq
\sum_{k\in \Z}A^{J(x)}P_k f(x)=\sum_{k\in \Z, k\ge K(x)}A^{J(x)}P_k f(x)+ \sum_{k\in \Z, k<K(x)}A^{J(x)}P_k f(x).
\endeq
For the latter part, i.e. the low frequency part, Bourgain's proof is already geometric, see Lemma 4.12 and Lemma 5.7 in \cite{Bo}. Hence the main task for us is to bound the former part, i.e. the high frequency part, by a geometric argument.

\begin{rem}
The estimate of the above high frequency part is done in Lemma 3.28 in \cite{Bo} via analytic methods. 
\end{rem}

We proceed with the estimate of the high frequency part: First we write
\begineq
\sum_{k\in \Z, k\ge K(x)}A^{J(x)}P_k f(x)=\sum_{l\in \N_0}A^{J(x)}P_{K(x)+l}f(x).
\endeq
Then by the triangle inequality, it suffices to prove that
\begineq
\|A^{J(x)}P_{K(x)+l}f(x)\|_2 \lesim 2^{-\mu l}\|f\|_2,
\endeq
for some $\mu>0$, with a constant being independent of $l\in \N$.\\

%\section{Replacing the $l^{\infty}$ sum by $l^2$ sume}
Notice that the above estimate is still of a maximal type, and we want to get rid of the linearization by replacing the $l^{\infty}$ norm by an $l^2$ norm. To do this, we need to introduce several notations. For $j, k\in \N$, define
\begineq
\Omega_{j, k}:=\{x\in \Omega | 2^{-j}\cdot \sup_{|t|<2^{-j}}\omega_x(t) \sim 2^{-k}\}.
\endeq 
For a real analytic vector field, either the integral curves are straight lines, or for each $j\in \N$, the complement of the set $\cup_{k}\Omega_{j,k}$ has measure zero. Hence it is no restriction to assume for each $j\in \N$ that
\begineq
\Omega=\bigcup_{k\in \N} \Omega_{j, k}.
\endeq
It is also clear that for a fixed $k\in \N$, the $\Omega_{j,k}$ for different $j$ are essentially disjoint. Hence for a fixed $x$, 
\begineq
|A^{J(x)}P_{K(x)+l}f(x)|\lesim \sup_{j\in \N} \left( \sum_{k\in \N}\left| A^j P_{k+l}f \right|^2 \mathbbm{1}_{\Omega_{j,k}} \right)^{1/2}.
\endeq
We replace the $\sup$ norm by the $l^2$ norm to obtain
\begineq\label{EE3.12}
|A^{J(x)}P_{K(x)+l}f(x)|\lesim \left(\sum_{j\in \N}\sum_{k\in \N}\left| A^j P_{k+l}f \right|^2 \mathbbm{1}_{\Omega_{j,k}} \right)^{1/2}.
\endeq
Taking the $L^2$ norm of \eqref{EE3.12}, we obtain
\begineq
\sum_{j\in \N}\sum_{k\in \N} \int \left| A^j P_{k+l}f \right|^2 \mathbbm{1}_{\Omega_{j,k}}.
\endeq 

%\vspace{4mm}

Recall that $R_{x, j}$ denotes the rectangle with center $x$, length $2^{-j}$ and width
\begineq
\delta(R_{x, j}):=2^{-j}\sup_{|t|<2^{-j}}\omega_x(t).
\endeq

%We plan to cover $\Omega_{j,k}$ by a collection of rectangles of the form $R_{x, j}$. Before doing this, let us first study some geometric properties of these rectangles.

\noindent In the following we will cover $\Omega_{j,k}$ by rectangles $\{R_m=R_{x_m, j}\}_{m\in \N}$ satisfying the following two conditions
\begineq
\begin{split}
& (i) \text{  } \text{  }\delta(R_m)=2^{-k};\\
& (ii) \text{  }\text{  }\text{the center of $R_m$ does not belong to } R_1\cup ... \cup R_{m-1}.
\end{split}
\endeq
%Moreover, we define
%\begineq
%\Omega'_{j,k}:=\bigcup_{x\in \Omega_{j,k}}2\cdot R_{x, {j}}.
%\endeq
Hence for fixed $j$ and $k$,
\begineq\label{EE3.18}
\int \left| A^j P_{k+l}f \right|^2 \mathbbm{1}_{\Omega_{j,k}}\lesim \sum_{m\in \N}\int \left| A^j P_{k+l}f \right|^2 \mathbbm{1}_{R_m}.
\endeq

Indeed, the above covering of $\Omega_{j,k}$ is a ``valid'' covering, i.e. a covering without much overlapping. To be precise, if we define
\begineq
\Omega'_{j,k}:=\bigcup_{x\in \Omega_{j,k}}2\cdot R_{x, {j}},
\endeq
then it has been proved by Bourgain in \cite{Bo} (see the following Lemma \ref{lemma3.6}) that
\begineq\label{AA3.20}
\|\sum_j \mathbbm{1}_{\Omega'_{j,k}}\|_{\infty}\lesim 1.
\endeq

In the following, when estimating the right hand side of \eqref{EE3.18}, we will need several other geometric properties like \eqref{AA3.20}. Hence we organize all them together in the next section.

%
%By Lemma \ref{lemma3.6}, we know that 
%\begineq
%\|\sum_{m\in \N} \mathbbm{1}_{2\cdot R_m}\|_{\infty}\lesim 1
%\endeq
%and therefore
%\begineq\label{FF3.18}
%\sum_{m\in \N}\mathbbm{1}_{2\cdot R_m}\lesim \mathbbm{1}_{\Omega'_{j,k}}.
%\endeq
%Indeed, Lemma \ref{lemma3.6} implies the following much stronger covering property:
%\begineq\label{AA3.20}
%\|\sum_j \mathbbm{1}_{\Omega'_{j,k}}\|_{\infty}\lesim 1.
%\endeq

\section{Geometric properties of rectangles}\label{section44}

In this section we collect several geometric lemmas that will play crucial roles in the forthcoming calculation. Most of these lemmas have already been proven by Bourgain \cite{Bo}. Here we still include them for the sake of completeness. Moreover, out of certain technical reasons, we will still need several variants of thess geometric lemmas, which are the following Lemma \ref{lemma3.4} and \ref{lemma3.7}.

\begin{lem}[Lemma 4.1 in \cite{Bo}]\label{lemma3.2}
Let $x'$ be in the rectangle $R_{x,j}$, then
\begineq
\delta(R_{x,j})\sim \delta(R_{x',j}),
\endeq
and  $R_{x,j}$ is contained in a multiple of  $R_{x',j}$ and vice versa.
\end{lem}

\begin{lem}[Lemma 4.6 in \cite{Bo}]\label{lemma3.3}
Assume 
\begineq
2\cdot R_{x,j}\cap 2\cdot R_{x',j'}\neq \emptyset,
\endeq
and $2^{-j'}\lesim 2^{-j}$. Then
\begineq
R_{x',j'}\subset 4\cdot R_{x,j}
\endeq
and 
\begineq
\frac{\delta(R_{x,j})}{2^{-j}}\gtrsim \frac{\delta(R_{x',j'})}{2^{-j'}}.
\endeq
Namely the larger rectangle has the larger eccentricity.
\end{lem}

\begin{lem}\label{lemma3.4}
Assume that 
\begineq
R_{x,j}\cap R_{x',j+j_0}\neq \emptyset,
\endeq
for some $j_0\in \N_0$. Then there exists a constant $a_0>1$ such that
\begineq\label{NN4.6}
\frac{\delta(R_{x,j})}{2^{-j}}\lesim (2^{j_0})^{a_0} \cdot \frac{\delta(R_{x',j+j_0})}{2^{-j-j_0}}.
\endeq
\end{lem}
\begin{rem}
Compared with Lemma \ref{lemma3.3}, this lemma says that the growth of the eccentricity of a rectangle with respect to its length can only be polynomial.
\end{rem}
\noindent {\bf Proof of Lemma \ref{lemma3.4}:} This follows easily from Lemma \ref{lemma3.2} and the following doubling estimate $(3.20)$ in \cite{Bo}
\begineq
\frac{\delta(R_{x,j})}{2^{-j}}\le C \cdot \frac{\delta(R_{x',j+1})}{2^{-j-1}},
\endeq
for some constant $C>0$. $\Box$\\

\begin{lem}[Lemma 4.7 in \cite{Bo}]\label{lemma3.6}
Let $\{R_{x_i, j_i}\}_{i\in \N_0}$ be a sequence of rectangles and $\delta>0$ such that 
\begineq
\begin{split}
& (i)  \text{ }\delta(R_{x_i, j_i})\sim \delta;\\
& (ii) \text{ } x_{i+1} \text{ does not belong to } R_{x_0,j_0}\cup ...\cup R_{x_i, j_i}, \forall i.
\end{split}
\endeq
Then
\begineq
\|\sum_{i\in \N_0}\mathbbm{1}_{2\cdot R_{x_i, j_i}}\|_{\infty}\lesim 1.
\endeq
\end{lem}

We will also need the following generalized version of the above lemma.
\begin{lem}\label{lemma3.7}
Under the same assumptions as in Lemma \ref{lemma3.6}, there exists a constant $b_0>0$ such that for any $N\in \N_0$, we have
\begineq
\|\sum_{i\in \N_0}\mathbbm{1}_{R^{p,q}_{x_i, j_i}}\|_{\infty}\lesim (p+q+1)^{b_0}, \forall p, q\in \N,
\endeq
where $R^{p,q}_{x_i, j_i}$ is obtained by dilating the length of $R_{x_i, j_i}$ to $p$ times, and the width to $q$ times.
\end{lem}
\begin{rem}
The above lemma says that when we enlarge the rectangles $R_{x_i,j_i}$, their overlapping can only be polynomially growing.
\end{rem}

\noindent {\bf Proof of Lemma \ref{lemma3.7}:} Denote
\begineq
\tilde{q}:=\max\{q, p^{a_0+1}\},
\endeq
where $a_0$ is the constant in \eqref{NN4.6}. For the rectangle $R^{p,q}_{x_i, j_i}$, we further enlarge it to be of width $\tilde{q}\cdot \delta$. Next, we will dilate the length to $\tilde{p}_i\cdot 2^{-j_i}$ such that 
\begineq
\sup_{|t|\le \tilde{p}_i\cdot 2^{-j_i}}\omega_{x_i}(t)\sim \frac{\tilde{q}}{\tilde{p}_i}\cdot \frac{\delta}{2^{-j_i}}.
\endeq
By Lemma \ref{lemma3.3}, it is not difficult to see that 
\begineq\label{NN4.13}
\tilde{p}_i\le \tilde{q}\lesim \tilde{p}_i^{a_0+1},
\endeq
uniformly in $i$.\\

Our goal now is to show that 
\begineq\label{NN4.14}
\|\sum_i \mathbbm{1}_{R^{\tilde{p}_i, \tilde{q}}_{x_i, j_i}}\|_{\infty}\lesim \tilde{q}^{b_0},
\endeq
for some $b_0$ to be determined later. Suppose that the $L^{\infty}$ norm on the left hand side of the above expression is attained at the point $O$. Moreover, let $\mathcal{R}_O$ denote the collection of rectangles containing the point $O$. W.l.o.g. we assume that
\begineq
\mathcal{R}_O=\{R^{\tilde{p}_i, \tilde{q}}_{x_i, j_i}\}_{0\le i\le N},
\endeq
for some $N\in \N_0$. Then \eqref{NN4.14} is equivalent to proving
\begineq\label{NN4.16}
N\lesim \tilde{q}^{b_0}.
\endeq
By the definition of the rectangles $R^{\tilde{p}_i, \tilde{q}}_{x_i, j_i}$ and Lemma \ref{lemma3.2}, it is not difficult to see that all the rectangles in $\mathcal{R}_O$ have comparable lengths. Indeed, up to a constant dilation factor, any of these rectangles is contained in another. Hence 
\begineq
\bigcup_{0\le i\le N}R^{\tilde{p}_i, \tilde{q}}_{x_i, j_i}\subset 4 R^{\tilde{p}_0, \tilde{q}}_{x_0, j_0}.
\endeq
Moreover, by the upper bound on $\tilde{p}_i$ in \eqref{NN4.13}, we can also obtain that 
\begineq
l(R_{x_i, j_i})\gtrsim l(R^{\tilde{p}_0, \tilde{q}}_{x_0, j_0})/\tilde{q},
\endeq
where for a rectangle $R$, $l(R)$ is used to denote its length. Hence by the assumption that the center $x_i$ of $R_{x_i, j_i}$ is not contained in 
\begineq
R_{x_0,j_0}\cup ...\cup R_{x_{i-1}, j_{i-1}}
\endeq
for all $i$, we obtain easily the estimate \eqref{NN4.16} for some constant $b_0$ depending only on $a_0$. So far we have finished the proof of Lemma \ref{lemma3.7}. $\Box$

\section{Estimate on each rectangle}

We proceed with the estimate of the right hand side of \eqref{EE3.18}. In this section, we will prove an estimate for fixed $j, k, l$ and for a given rectangle $R_m$. Here $R_m$ is a rectangle of length $2^{-j}$ and width $2^{-k}$. \\

Recall that we have assumed that the vector field points in the cone $\Gamma_0$, which is the two-ended cone forming an angle less than $\pi/10$ with the horizontal axis. If we denote by $P_{\Gamma_0}$ the frequency projection operator on the cone $\Gamma_0$, it is not difficult to see that 
\begineq
A^j P_{k+l}P_{\Gamma_0}f\equiv 0.
\endeq
Hence in the following we will only be concerned with the frequency in the cone $\Gamma_0^c$. Moreover, for the sake of simplicity, we will always identify $P_{k+l}$ with $P_{k+l}P_{\Gamma_0^c}$.

\subsection{Time-frequency decomposition of the function $P_{k+l}f$}

Most of the content in this subsection is taken from Bateman \cite{Ba2}. Here we will make some modifications as we will be dealing with all frequency annuli instead of one single annulus.\\

{\bf Frequency decomposition.} For the fixed $j,k$ and $l$, we will denote
\begineq
\theta:=k+l-j.
\endeq 
Let $\mathcal{D}_{\theta}$ be the collection of the dyadic intervals of length $2^{-\theta}$ contained in $[-20,20]$. Fix a smooth positive function $\beta:\R\to \R$ s.t. 
\begineq
\beta(x)=1, \forall |x|\le 1; \beta(x)=0, \forall |x|\ge 2.
\endeq
Also choose $\beta$ such that $\sqrt{\beta}$ is a smooth function. For each $\omega\in \mathcal{D}_{\theta}$, define
\begineq
\beta_{\omega}(x)= \beta(2^{\theta}(x-c_{\omega})),
\endeq
where $c_{\omega}$ denotes the center of the interval $\omega$. Define
\begineq
\beta_{\theta}(x)=\sum_{\omega\in \mathcal{D}_{\theta}} \beta_{\omega}(x).
\endeq
Notice that 
\begineq
\beta_{\theta}(x+2^{-\theta})=\beta_{\theta}(x), \forall x\in [-20, 20-2^{-\theta}].
\endeq
Define 
\begineq
\gamma_{\theta}=\frac{1}{20}\int_{-10}^{10} \beta_{\theta}(x+t)dt.
\endeq
Because of the above periodicity, we know that $\gamma_{\theta}$ is constant for $x\in [-10,10]$, independently of $\theta$. Say $\gamma_{\theta}(x)=\gamma>0$, hence 
\begineq
\frac{1}{\gamma}\cdot \gamma_{\theta}(x) \mathds{1}_{[-10,10]}(x)=\mathds{1}_{[-10,10]}(x).
\endeq
Define another multiplier $\tilde{\beta}:\R\to \R$ with support in $[\frac{1}{2}, \frac{5}{2}]$ and $\tilde{\beta}(x)=1$ for $x\in [1,2]$. We define the corresponding multipliers on $\R^2$:
\begin{align*}
& \hat{m}_{k+l,\omega}(\xi, \eta)=\tilde{\beta}(2^{-k-l}\eta)\beta_{\omega}(\frac{\xi}{\eta})\\
& \hat{m}_{k+l,\theta,t}(\xi, \eta)=\tilde{\beta}(2^{-k-l}\eta)\beta_{\theta}(t+\frac{\xi}{\eta})\\
& \hat{m}_{k+l,\theta}(\xi, \eta)=\tilde{\beta}(2^{-k-l}\eta)\gamma_{\theta}(\frac{\xi}{\eta})
\end{align*}

%The following picture illustrates the decomposition in frequency.

%\begin{minipage}[b]{10.2cm}
%\includegraphics[width=10cm]{frequencydecomposition.png}
%\end{minipage}

\noindent Using the above multipliers, we obtain
\begineq
A^j P_{k+l}f=A^j (m_{k+l,\theta}* f)=\frac{1}{20}\int_{-10}^{10}A^j (m_{k+l, \theta, t}* f)dt.
\endeq
Hence it suffices to prove a uniform bound on $t\in [-10, 10]$. W.l.o.g. we will just consider the case $t=0$, which is 
\begineq\label{AA3.40}
A^j (m_{k+l, \theta, 0}* f)=\sum_{\omega\in \mathcal{D}_{\theta}}A^j (m_{k+l, \omega}* f).
\endeq

{\bf Space (Time) decomposition:} For $\omega\in \mathcal{D}_{\theta}$, let $\mathcal{U}_{k+l,\omega}$ be a partition of $\R^2$ by rectangles of width $2^{-k-l}$ and length $2^{-j}$, whose long side have slope $-c_{\omega}$ with $c_{\omega}$ denoting the center of the interval $\omega$. If $s\in \mathcal{U}_{k+l,\omega}$, we will write $\omega_s:=\omega$.

\begin{defi}
For a rectangle $R\subset\R^2$ of slope less than one, with $l(R)$ its length, $w(R)$ its width, we define its  uncertainty interval $EX(R)$ to be the interval of width $w(R)/l(R)$ and centered at $slope(R)$.
\end{defi}

\begin{rem}
For a tile $s\in \mathcal{U}_{k+l, \omega}$, we have that $EX(s)=-\omega$.
\end{rem}

%Let
%\begin{align*}
%& \mathcal{U}_{k,l}:= \cup_{\omega\in \mathcal{D}_{k-l}}\mathcal{U}_{k,\omega},
%\end{align*}

An element of $\mathcal{U}_{k+l, \omega}$ for some $\omega\in \mathcal{D}_{\theta}$ is called a ``tile''. Define $\varphi_{k+l,\omega}$ such that
\begineq
|\hat{\varphi}_{k+l,\omega}|^2=\hat{m}_{k+l,\omega},
\endeq
then $\varphi_{k+l,\omega}$ is smooth by our assumption on $\beta$ mentioned above.\\

For a tile $s\in \mathcal{U}_{k+l,\omega}$, define 
\begineq\label{wavelet}
\varphi_s(p):= \sqrt{|s|} \varphi_{k+l,\omega} (p-c_s),
\endeq
where $c_s$ is the center of $s$. Notice that 
\begineq
\|\varphi_s\|_2^2= \int_{\R^2} |s| \varphi_{k+l,\omega}^2=|s|\int_{\R^2} \hat{m}_{k+l,\omega}=1,
\endeq
i.e. $\varphi_s$ is $L^2$ normalized.\\

The construction of the tiles above by the uncertainty principle is to localize functions further in space, for this purpose we need
\begin{lem}(\cite{Ba2})
Under the above notations, for the frequency localised function $f*m_{k+l, \omega}$, we have 
\begineq
f*m_{k+l,\omega}(x)=\lim_{N\to \infty} \frac{1}{4 N^2} \int_{[-N,N]^2} \sum_{s\in \mathcal{U}_{k+l,\omega}} \langle f,\varphi_s(p+\cdot)\rangle \varphi_s(p+x)dp.
\endeq
\end{lem}
\noindent The above lemma allows us to pass the expression in \eqref{AA3.40} to the model sum
\begin{align*}
\sum_{\omega\in \mathcal{D}_{\theta}}\sum_{s\in \mathcal{U}_{k+l, \omega}}\langle f, \varphi_s\rangle A^j \varphi_s.
\end{align*}

\begin{defi}
For a unit vector $v=(v_1, v_2)\in \R^2$, for an interval $I\subset \R$, we say that $v\in I$ if the slope of the line that is perpendicular to $v$ lies in $I$.
%\begineq
%-\frac{v_1}{v_2}\in I.
%\endeq
\end{defi}
\begin{lem}\label{lemma5}(\cite{LL1}, \cite{Ba2})
we have that $A^j \varphi_s(x)=0$ unless $v(x)\in \omega_{s}$.
\end{lem}
To have this lemma is the main reason of replacing the strict cut-off by a smooth cut-off in Section \ref{section2}. The proof of Lemma \ref{lemma5} is simply by the Plancherel theorem. From this lemma we know that in order for the output $A^j \varphi_s$ not to vanish at a point $x$, the vector field at $x$ has to point roughly to the direction of the long side of the tile $s$.

\subsection{Estimate on each rectangle by ignoring the tails of the wavelet functions}

After the above preparation, we turn to the estimate of each term on the right hand side of \eqref{EE3.18}, which is
\begineq\label{AA3.47}
\int \left| A^j P_{k+l}f \right|^2 \mathbbm{1}_{R_m} = \int \left| \sum_{\omega\in \mathcal{D}_{\theta}}\sum_{s\in \mathcal{U}_{k+l, \omega}}\langle f, \varphi_s\rangle A^j \varphi_s \right|^2 \mathbbm{1}_{R_m}.
\endeq
By Lemma \ref{lemma5}, we have that the right hand side of \eqref{AA3.47} is equal to
\begineq\label{AA3.48}
\sum_{\omega\in \mathcal{D}_{\theta}} \int \left| \sum_{s\in \mathcal{U}_{k+l, \omega}}\langle f, \varphi_s\rangle A^j \varphi_s \right|^2 \mathbbm{1}_{R_m}.
\endeq

In this subsection, we will only show the ideas of how to bound the above term, or in another word, we will ignore the tails of the wavelet functions and the function $\alpha$ in the definition of $A^j$ in \eqref{AA2.1}, and always assume that they have compact support in both space and frequency.\\

Under the above simplification, the expression in \eqref{AA3.48} becomes
\begineq\label{AA3.49}
\sum_{\omega\in \mathcal{D}_{\theta}} \sum_{s\in \mathcal{U}_{k+l, \omega}} |\langle f, \varphi_s\rangle|^2 \int \left|  A^j \varphi_s \right|^2 \mathbbm{1}_{R_m}.
\endeq

\noindent Take a point $x\in R_m$, for a tile $s\in \mathcal{U}_{k+l, \omega}$ for some $\omega\in \mathcal{D}_{\theta}$, we observe that in order for $A^j \varphi_s(x)$ not to vanish, we must have $\omega\subset 3\cdot EX(R_m)$ as by Lemma \ref{lemma3.2} we know that $v(x)\in 2\cdot EX(R_m)$ for any $x\in R_m$. This, together with the fact that both $R_m$ and $s$ have length $2^{-j}$, implies that
\begineq
s\subset 4\cdot R_m,
\endeq
for those tiles $s$ such that $A^j \varphi_s$ is not identically zero.

\begin{claim}\label{claim3.14}
There exists $\mu>0$ such that 
\begineq
\int \left|  A^j \varphi_s \right|^2 \mathbbm{1}_{R_m}\lesim 2^{-\mu l},
\endeq
with the constant being independent of $s$.
\end{claim}

By the above claim, the expression in \eqref{AA3.49} can be further bounded by
\begineq
\sum_{\omega\in \mathcal{D}_{\theta}} \sum_{s\in \mathcal{U}_{k+l, \omega}, s\subset 4\cdot R_m} 2^{-\mu l} \cdot  |\langle f, \varphi_s\rangle|^2 \lesim 2^{-\mu l} \|\mathbbm{1}_{R_m}\cdot P_{k+l}f\|_2^2.
\endeq
The next step is to sum over $m, j$ and $k$:
\begineq
\begin{split}
\sum_{j,k}\sum_{m}\|\mathbbm{1}_{R_m}\cdot P_{k+l}f\|_2^2  &  \lesim 2^{-\mu l} \sum_{j,k}\|\mathbbm{1}_{\Omega'_{j,k}}\cdot P_{k+l}f\|_2^2\\
& \lesim 2^{-\mu l}\sum_{k}\|P_{k+l}f\|_2^2 \lesim 2^{-\mu l}\|f\|_2^2,
\end{split}
\endeq
where we have used the disjointness property \eqref{AA3.20}. Hence for the model problem, what remains is\\

{\bf ``Proof'' of Claim \ref{claim3.14}:} We can w.l.o.g. assume that there exists a point $x_0\in s$ such that
\begineq
v(x_0)\in \omega_s,
\endeq
as otherwise $A^j \varphi_s$ will be identically zero. By a proper translation and rotation, we can assume that $x_0=(0,0)$ and $v(x_0)=(1, 0)$. 

Now we look at the directions of the vector field at the points lying on the line segment 
\begineq
\{(x_1, x_2): x_2=0\}\cap s.
\endeq
By the assumption on the rectangle $R_m$ we know that
\begineq
\sup_{|t|\le 2^{-j}} w_{x_0}(t)=\sup_{|t|\le 2^{-j}} |\det[v(x_0+t v(x_0)), v(x_0)]|\sim 2^{-k+j}.
\endeq
Notice that $|\omega_s|=2^{-k-l+j}$, hence in order for $A^j \varphi_s$ not to vanish at a point $x\in s\cap \{(x_1, x_2): x_2=0\}$, we must have
\begineq
|\det[v(x), v(x_0)]|=w_{x_0}(x\cdot v(x_0))\lesim 2^{-k-l+j}.
\endeq
By taking $\tau=2^{-l}$ in the condition \eqref{AA1.3} we obtain
\begineq
\left| { t\in [-2^{-j}, 2^{-j}]: w_{x_0}(t)< 2^{-l} \sup_{|t|\le 2^{-j}} w_{x_0}(t) } \right|\le C_0 2^{-c_0 l} \cdot 2^{-j},
\endeq
which further implies that 
\begineq
\left| \{(x_1, 0)\in s: A^j \varphi_s(x_1, 0)\neq 0)\} \right| \lesim 2^{-c_0 l} \cdot 2^{-j}.
\endeq
So far we have proved that on one line segment, the non-vanishing output $A^j \varphi_s$ has relatively small measure. In the next, we want to show that this indeed holds true for all the points in the tile $s$, namely
\begineq\label{AA3.60}
|\{x\in s: A^j \varphi_s(x)\neq 0\}|\lesim 2^{-c_0 l}|s|.
\endeq
This, combined with the trivial estimate
\begineq
\|A^j \varphi_s\|_{\infty}\lesim |s|^{-1/2}, 
\endeq
concludes the proof of Claim \ref{claim3.14}.\\

We turn to the proof of \eqref{AA3.60}: For $|x_2|\le 2^{-k-l+2}$, consider the line segment
\begineq
L_{x_2}:= \{(0, x_2)+t\cdot v(0, x_2): |t|\le 2^{-j+2}\}.
\endeq
First by the $C^1$ assumption on the vector field, we know that 
\begineq
v(0, x_2)\in 2\cdot \omega_s, \forall |x_2|\le 2^{-k-l+2}.
\endeq
Then by the same argument as before, we obtain that 
\begineq
|\{x\in L_{x_2}: A^j \varphi_s(x)\neq 0\}| \lesim 2^{-c_0 l}\cdot 2^{-j},
\endeq
for each $|x_2|\le 2^{-k-l+2}$. Hence by Fubini's theorem (which can be applied due to the $C^1$ assumption on the vector field), we obtain
\begineq
\begin{split}
& |\{x\in s: A^j \varphi_s(x)\neq 0\}|\\
&  =\int_{-2^{-k-l+2}}^{2^{-k-l+2}} |\{x\in L_{x_2}: A^j \varphi_s(x)\neq 0\}| dx_2 \lesim 2^{-c_0 l}\cdot 2^{-k-l-j}.
\end{split}
\endeq
Hence we have finished the proof of \eqref{AA3.60}.

\subsection{The full estimate on each rectangle}

In this subsection, we will make the above heuristic argument rigorous, i.e. we will also take care of the tails of the wavelet functions. For fixed $j, k, l$ and $m$, we want to bound the following
\begineq
\int \left| A^j P_{k+l}f \right|^2 \mathbbm{1}_{R_m} = \sum_{\omega\in \mathcal{D}_{\theta}} \int \left| \sum_{s\in \mathcal{U}_{k+l, \omega}}\langle f, \varphi_s\rangle A^j \varphi_s \right|^2 \mathbbm{1}_{R_m}.
\endeq
For $p, q\in \Z$, we denote by $\vec{R}_m^{p,q}$ the translation of the rectangle $R_m$ by $(p, q)$ units, i.e.
\begineq
\vec{R}_m^{p,q}=R_{m}+p\cdot 2^{-j} v(x_m)+ q\cdot 2^{-k} v^{\perp}(x_m),
\endeq
where $x_m$ denotes the center of $R_m$ and $v(x_m)$ is the value of the vector field at the point $x_m$ which is parallel to the long side of $R_m$.

Hence for one fixed $\omega\in \mathcal{\theta}$, we have
\begineq
\int \left| \sum_{s\in \mathcal{U}_{k+l, \omega}}\langle f, \varphi_s\rangle A^j \varphi_s \right|^2 \mathbbm{1}_{R_m}=\int \left| \sum_{p,q\in \Z}\sum_{s\in \mathcal{U}_{k+l, \omega}, s\subset \vec{R}_{m}^{p,q}}\langle f, \varphi_s\rangle A^j \varphi_s \right|^2 \mathbbm{1}_{R_m}.
\endeq
By Minkowski's inequality, the right hand side of the above display can be bounded by
\begineq\label{AA3.69}
\left(\sum_{p,q\in \Z} \left(\int \left| \sum_{s\in \mathcal{U}_{k+l, \omega}, s\subset \vec{R}_{m}^{p,q}}\langle f, \varphi_s\rangle A^j \varphi_s \right|^2 \mathbbm{1}_{R_m}\right)^{1/2}\right)^2.
\endeq
\begin{lem}\label{lemma3.15}
For any large $M\in \N_0$, there exists a constant $C_M$ such that
\begineq
\int \left| \sum_{s\in \mathcal{U}_{k+l, \omega}, s\subset \vec{R}_{m}^{p,q}}\langle f, \varphi_s\rangle A^j \varphi_s \right|^2 \mathbbm{1}_{R_m} \lesim \frac{C_M\cdot 2^{-\mu l}}{(|p|+|q|+1)^M}\sum_{s\in \mathcal{U}_{k+l, \omega}, s\subset \vec{R}_{m}^{p,q}} |\langle f, \varphi_s\rangle|^2,
\endeq
where $\mu$ is the same as the one in Claim \ref{claim3.14}.
\end{lem}
We substitute the estimate in Lemma \ref{lemma3.15} into \eqref{AA3.69} to obtain
\begineq
\begin{split}
& \left( \sum_{p,q\in \Z} \left( \frac{2^{-\mu l}}{(|p|+|q|+1)^M}\sum_{s\in \mathcal{U}_{k+l, \omega}, s\subset \vec{R}_{m}^{p,q}} |\langle f, \varphi_s\rangle|^2 \right)^{1/2} \right)^2\\
& \lesim \sum_{p,q\in \Z} \frac{2^{-\mu l}}{(|p|+|q|+1)^M}\sum_{s\in \mathcal{U}_{k+l, \omega}, s\subset \vec{R}_{m}^{p,q}} |\langle f, \varphi_s\rangle|^2,
\end{split}
\endeq
where the exact value of $M$ might vary from line to line. Hence we have obtained
\begineq\label{AA3.72}
\int \left| A^j P_{k+l}f \right|^2 \mathbbm{1}_{R_m} \lesim \sum_{\omega\in \mathcal{D}_{\theta}}\sum_{p,q\in \Z} \frac{2^{-\mu l}}{(|p|+|q|+1)^M}\sum_{s\in \mathcal{U}_{k+l, \omega}, s\subset \vec{R}_{m}^{p,q}} |\langle f, \varphi_s\rangle|^2,
\endeq
which is the estimate on one single rectangle that we are aiming at.\\

{\bf Proof of Lemma \ref{lemma3.15}:} We will only consider the case $p=q=0$. The decay in $p$ and $q$ in the other case will simply follow from the non-stationary phase method. Hence what we need to prove becomes
\begineq
\int \left| \sum_{s\in \mathcal{U}_{k+l, \omega}, s\subset R_{m}}\langle f, \varphi_s\rangle A^j \varphi_s \right|^2 \mathbbm{1}_{R_m} \lesim 2^{-\mu l}\sum_{s\in \mathcal{U}_{k+l, \omega}, s\subset R_{m}} |\langle f, \varphi_s\rangle|^2.
\endeq
Recall that each tile $s$ has width $2^{-k-l}$, however the rectangle $R_m$ has width $2^{-k}$. This suggests that we should do a further partition of $R_m$ into smaller rectangles which will be of the same dimension as $s$.

We enumerate the tiles $s\subset R_m$ from above to below by $s_1, s_2, ..., s_{m'}...$, where $m'\lesim 2^l$. Notice that
\begineq
R_m\subset \bigcup_{m'} 2\cdot s_{m'}.
\endeq
Hence 
\begineq
\int \left| \sum_{s\in \mathcal{U}_{k+l, \omega}, s\subset R_{m}}\langle f, \varphi_s\rangle A^j \varphi_s \right|^2 \mathbbm{1}_{R_m} \lesim \sum_{m'}\int \left| \sum_{m^{''}}\langle f, \varphi_{s_{m^{''}}}\rangle A^j \varphi_{s_{m^{''}}} \right|^2 \mathbbm{1}_{2\cdot s_{m'}}.
\endeq
By the Cauchy-Schwartz inequality, we bound the right hand side of the above expression by
\begineq\label{AA3.76}
\sum_{m'}\int \sum_{m^{''}}  \left| \langle f, \varphi_{s_{m^{''}}}\rangle A^j \varphi_{s_{m^{''}}} \right|^2 (|m'-m^{''}|+1)^M \mathbbm{1}_{2\cdot s_{m'}},
\endeq
for some large constant $M$. By the same argument as in the proof of Claim \ref{claim3.14}, we obtain
\begineq
|\{x\in 2\cdot s_{m^{'}}: A^j \varphi_{s_{m^{''}}}\neq 0\}|\lesim 2^{-c_0 l} |s_{m'}|.
\endeq
This, together with the trivial bound
\begineq
\|A^j \varphi_{s_{m^{''}}}\|_{L^{\infty}(s_{m'})}\lesim \frac{|s_{m^{''}}|^{1/2}}{(|m'-m^{''}|+1)^{2M}},
\endeq
implies that 
\begineq
\int \left|\langle f, \varphi_{s_{m^{''}}}\rangle A^j \varphi_{s_{m^{''}}}\right|^2 \mathbbm{1}_{2\cdot s_{m'}} \lesim \frac{2^{-c_0 l}}{(|m'-m^{''}|+1)^{2M}} |\langle f, \varphi_{s_{m^{''}}}\rangle|^2.
\endeq
We substitute the above estimate into \eqref{AA3.76} to obtain
\begineq
 \sum_{m'}  \sum_{m^{''}} \frac{2^{-c_0 l}}{(|m'-m^{''}|+1)^M} |\langle f, \varphi_{s_{m^{''}}}\rangle|^2 \lesim \sum_{m^{''}} |\langle f, \varphi_{s_{m^{''}}}\rangle|^2.
\endeq
So far we have finished the proof of Lemma \ref{lemma3.15}, hence the estimate on each rectangle, which is \eqref{AA3.72}.

\section{Organizing all rectangles to finish the proof}

In this section, we will organize the estimates on all the rectangles together, i.e. to finish the proof of the following estimate
\begineq\label{AA3.81}
\sum_{j,k} \sum_{m}\int \left| A^j P_{k+l}f \right|^2 \mathbbm{1}_{R_m} \lesim 2^{-\mu l}\|f\|_2^2,
\endeq
for some $\mu>0$. To do this, we substitute the estimate \eqref{AA3.72} into the left hand side of the above expression to obtain
\begineq
\sum_{k,j}\sum_m \sum_{\omega\in \mathcal{D}_{\theta}}\sum_{p,q\in \Z} \frac{2^{-\mu l}}{(|p|+|q|+1)^M}\sum_{s\in \mathcal{U}_{k+l, \omega}, s\subset \vec{R}_{m}^{p,q}} |\langle f, \varphi_s\rangle|^2,
\endeq
where we are still using the notation $\theta=k+l-j$. Hence it suffices to prove that for fixed $p, q\in \Z$, we have
\begineq\label{AA3.83}
\sum_{k,j}\sum_m \sum_{\omega\in \mathcal{D}_{\theta}}\sum_{s\in \mathcal{U}_{k+l, \omega}, s\subset \vec{R}_{m}^{p,q}} |\langle f, \varphi_s\rangle|^2 \lesim (|p|+|q|+1)^{b_0} \|f\|_2^2,
\endeq
where $b_0$ is the constant in Lemma \ref{lemma3.7}.\\

{\bf Proof of the estimate \eqref{AA3.83}:} We first fix $k$. For the case $p=q=0$, for two tiles $s'$ and $s^{''}$ in the following collection 
\begineq
\bigcup_j \bigcup_{\omega\in \mathcal{D}_{\theta}}\bigcup_m \{s: s\in \mathcal{U}_{k+l, \omega}, s\subset R_m\},
\endeq
we either have 
\begineq
\omega_{s'}\cap \omega_{s^{''}}=\emptyset,
\endeq
or 
\begineq\label{AA3.86}
s'\cap s^{''} =\emptyset.
\endeq
Hence by the (almost) orthogonality, we obtain that 
\begineq
\sum_{j}\sum_m \sum_{\omega\in \mathcal{D}_{\theta}}\sum_{s\in \mathcal{U}_{k+l, \omega}, s\subset R_{m}} |\langle f, \varphi_s\rangle|^2 \lesim \|P_{k+l}f\|_2^2.
\endeq
By summing over $k$, we get the desired estimate \eqref{AA3.83} for the case $p=q=0$.\\

For the general $p,q\in \Z$, we no longer have \eqref{AA3.86} due to the simple fact that for two disjoint rectangles (of different dimensions), they might intersect after being translated by $(p,q)$ units separately. Fortunately, Lemma \ref{lemma3.7} says that the intersection caused by translation can only grow polynomially in $p$ and $q$.

Hence by essentially the same idea as above and by losing a factor of $(|p|+|q|+1)^{b_0}$, we obtain
\begineq
\sum_{j}\sum_m \sum_{\omega\in \mathcal{D}_{\theta}}\sum_{s\in \mathcal{U}_{k+l, \omega}, s\subset \vec{R}_{m}^{p,q}} |\langle f, \varphi_s\rangle|^2 \lesim (|p|+|q|+1)^{b_0} \|P_{k+l}f\|_2^2.
\endeq
Summing over $k$, we get the estimate \eqref{AA3.83}. So far we have finished the proof of \eqref{AA3.81}, hence the geometric proof of Theorem \ref{theorem1.1}.

%The crucial lemma
%\begin{lem}\label{lemma3.1}
%Under the above notations, we have
%\begineq\label{EE3.19}
% \|\mathbbm{1}_{R_m}\cdot A^j P_{k+l}f  \|_2 \lesim 2^{-\gamma l}\|P_{k+l}f \cdot \left(\mathbbm{1}_{2\cdot R_m}* \psi_{k+l}\right)\| _2,
%\endeq
%for some $\gamma>0$, which is independent of $k, j$ and $m$.
%\end{lem}
%
%We substitute the estimate in \eqref{EE3.19} into \eqref{EE3.18} to obtain
%\begineq
%\begin{split}
%\int \left| A^j P_{k+l}f \right|^2 \mathbbm{1}_{\Omega_{j,k}} &\lesim \sum_{m\in \N} 2^{-2\gamma l}\int \left| P_{k+l}f\right|^2 \left(\mathbbm{1}_{2\cdot R_m}* \psi_{k+l}\right)^2\\
%										&\lesim \sum_{m\in \N} 2^{-2\gamma l}\int \left| P_{k+l}f\right|^2 \left(\mathbbm{1}_{2\cdot R_m}* \psi_{k+l}\right).
%\end{split}
%\endeq
%By using \eqref{FF3.18}, we bound the last expression by
%\begineq
%2^{-2\gamma l}\int \left| P_{k+l}f\right|^2 \left(\mathbbm{1}_{\Omega'_{j,k}}* \psi_{k+l}\right).
%\endeq
%In the end, we just need to sum over $j$ and $k$: observe that for any $k$, we have
%\begineq
%\|\sum_{j} \mathbbm{1}_{\Omega'_{j,k}}\|_{\infty}\lesim 1,
%\endeq
%with some universal constant. Hence
%\begineq
%\begin{split}
%& \sum_{j,k}2^{-2\gamma l}\int \left| P_{k+l}f\right|^2 \left(\mathbbm{1}_{\Omega'_{j,k}}* \psi_{k+l}\right)\\
%& \lesim \sum_{k\in \N}2^{-2\gamma l}\int \left| P_{k+l}f\right|^2 \lesim 2^{-2\gamma l} \|f\|_2.
%\end{split}
%\endeq
%
%\noindent {\bf Proof of Lemma \ref{lemma3.1}:}

\vspace{1cm}

Shaoming Guo, Institute of Mathematics, University of Bonn\\
\indent Address: Endenicher Allee 60, 53115, Bonn\\
\indent Email: shaoming@math.uni-bonn.de

\end{document}